\documentclass[10pt]{article}
\usepackage{amssymb,amsmath,amsthm,hyperref}

\usepackage{geometry}
\newtheorem{prop}{Proposition}

\newtheorem{theorem}[prop]{Theorem}

\newtheorem{conjecture}[prop]{Conjecture}

\newtheorem{question}[prop]{Question}

\theoremstyle{definition}
\newtheorem{defn}[prop]{Definition}

\newcommand{\B}{\ensuremath{\mathbb{B}}} 
\newcommand{\Bexp}{\ensuremath{\B_{\exp}}}

\newcommand{\Z}{\ensuremath{\mathbb{Z}}}
\newcommand{\Q}{\ensuremath{\mathbb{Q}}}

\newcommand{\C}{\ensuremath{\mathbb{C}}}
\newcommand{\Cexp}{\ensuremath{\mathbb{C}_{\mathrm{exp}}}}

\newcommand{\M}{\mathcal{M}}

\newcommand{\ga}{\ensuremath{\mathbb{G}_\mathrm{a}}}   
\newcommand{\gm}{\ensuremath{\mathbb{G}_\mathrm{m}}}  

\newcommand{\tuple}[1]{\ensuremath{\langle #1 \rangle}}
\newcommand{\class}[2]{\ensuremath{\left\{ #1 \,\left|\, #2 \right.\right\}}}

\newcommand{\leteq}{\mathrel{\mathop:}=}
\newcommand{\iso}{\cong}

\DeclareMathOperator{\td}{td}  
\DeclareMathOperator{\ldim}{ldim}  

\title{Around Zilber's quasiminimality conjecture}
\author{Jonathan Kirby}

\begin{document}
\maketitle
This is an extended abstract for a survey talk given in Oberwolfach on 1st December 2022, slightly updated in June 2023.

\medskip

About 25 years ago, Zilber stated:

\begin{conjecture}\cite{Zilber97}\label{QM conj}
The complex field with the exponential function, $\tuple{\C;+,\cdot,\exp}$, is quasiminimal (QM): every definable subset is countable or co-countable. (Definable here means definable with parameters.)
\end{conjecture}

This conjecture has sparked a lot of activity over that time. For example, Zilber's part of the Zilber-Pink conjecture and the related work on functional transcendence came out of his early work towards the quasiminimality conjecture. Recently there has been significant progress towards proving the conjecture itself.

In this talk I surveyed some of the work around the conjecture, including the recent result of Gallinaro and myself that the complex field equipped with complex power functions is quasiminimal.
\begin{theorem}\cite{GallinaroKirby23}\label{powers theorem} \ 
For $\lambda \in \C$, let $\Gamma_\lambda = \class{(\exp(z), \exp(\lambda z))}{z \in \C}$, the graph of the multivalued map $w \mapsto w^\lambda$. Then the structure
\[\tuple{\C;+,\cdot,-,0,1,(\Gamma_\lambda)_{\lambda \in \C}}\]
 is quasiminimal.
\end{theorem}

\section{Quasiminimality}
The complex field $\C_{\mathrm{field}}$ is well-known to be \emph{minimal}: every definable subset of $\C$ is finite or co-finite, and indeed \emph{strongly minimal}: the same is true of every model of its first-order theory, which is $\mathrm{ACF}_0$.

It is well-known that strongly minimal theories are uncountably categorical, (model-theoretic) algebraic closure is a pregeometry, and all models are prime over a basis of the pregeometry.

Quasiminimality as defined in Conjecture~\ref{QM conj} corresponds to minimality. There is a stronger property:
\begin{defn}
An uncountable structure $\M$ is \emph{quasiminimal in the sense of automorphisms} (QM$_{\text{Aut}}$) if for all reducts $\M_0$ of $\M$ to countable languages, and all countable subsets $A$ of (the underlying set of) $\M_0$, there is a co-countable orbit of $\mathrm{Aut}(\M_0/A)$.
\end{defn}
This still corresponds to minimality, not strong minimality. There is a yet stronger property, called Quasiminimal Excellence (QME), defined by Zilber \cite{Zilber05qmec}, improved in \cite{OQMEC} and \cite{Baldwin_Categoricity} and then substantially simplified in \cite{BHHKK14}. I refer the reader to the last of these for the simplest definition. Putting the work of these papers together we get
\begin{theorem}
If $\M$ is uncountable and QME then its $L_{\omega,\omega_1}(Q)$-theory is uncountably categorical, the ``countable closure'' operator is a pregeometry, and models of that theory are prime over bases of this pregeometry.
\end{theorem}

So we have QME $\implies$ QM$_\text{Aut}$ $\implies$ QM. A natural question was asked during the Oberwolfach meeting:
\begin{question}
Suppose $|\M| \geqslant \aleph_2$ and $\M$ is QM. Must it be QME?
\end{question}
A previously known counterexample for $|\M| = \aleph_1$ is given by the dense linear order $\omega_1 \times_{\text{lex}} \Q_<$, which is QM$_\text{Aut}$ but not QME, and indeed the ``countable closure'' operator is not a pregeometry. This structure is approximated by countable substructures in a fundamentally different way, via the linear order $\omega_1$, whereas in the excellent case one approximates by the directed partial order of all countable subsets of some uncountable basis.
Work of Pillay and Tanovic \cite{PillayTanovic} may be relevant here.

Subsequently it turned out that the answer to the question is no, there are examples of QM structures of arbitrarily large cardinality, for which countable closure is a pregeometry, but which are not QM$_\text{Aut}$. There also seem to be examples which are QM$_\text{Aut}$ but not QME. A revised question would be to understand what different behaviours are possible.

\section{Variant conjectures}
In Zilber's conjecture, one can replace the complex exponential function with a Weierstrass $\wp$-function, or equivalently by the exponential map of an elliptic curve. More generally, one could consider the exponential map of an abelian variety, or a semiabelian variety, or indeed any commutative complex algebraic group.

One could even consider the exponential maps of all such groups together in one structure. All these conjectures seem fairly equally plausible and equally difficult. The Ax--Schanuel theorem \cite{Ax71,Ax72a} gives a strong structural result in all these cases.

While we now have an Ax--Schanuel theorem for the modular $j$-function \cite{PilaTsimerman16}, the domain of that is the upper half-plane so $\tuple{\C;+,\cdot,j}$ and similar examples are far from quasiminimal.

Koiran asked if the expansion of the complex field by all unary entire complex functions is QM, and Wilkie 
noted that we have no counterexample. This appears to be a much harder conjecture than the group case.

The work of Fatou and Bieberbach from the 1920s mentioned in \cite{Zilber97} gives examples of holomorphic $f:\C^2 \to \C^2$ whose image is open but not dense, and they cannot be definable in a QM structure.

\section{Some results}

It is immediate that $\C_\Z \leteq \tuple{\C;+,\cdot,\Z}$, and even $\C_{\Z \mathrm{IP}} \leteq \tuple{\C;+,\cdot,\Z,(z,n) \mapsto z^n}$, the complex field with a predicate for $\Z$ and the function of raising complex numbers to integer powers, are QM$_\text{Aut}$, because they have the same automorphisms as the pure field, and indeed they are QME.

Zilber~\cite{Zilber02tgfd} gave a theory of a generic function on a field, which turns out to be the theory of an ultrapower of polynomials of degree tending to infinity~\cite{Koiran01}. Any such entire function is QM, for topological reasons. Wilkie~\cite{Wilkie05} used power series with sparse rational coefficients to construct entire functions, which he called Liouville functions after Liouville's construction of transcendental numbers, and partly showed they are generic in Zilber's sense, before Koiran \cite{Koiran03} finished the proof.

Towards Conjecture~\ref{QM conj}, Wilkie made progress around an Oberwolfach meeting in 2003, which led to an announcement of the proof of quasiminimality of raising to the power $i$ in 2008, although that is still unpublished. During a long visit to Oxford in 2004, Macintyre was working both on ideas to prove the conjecture and to provide a counterexample.

Boxall~\cite{Boxall20} showed that certain existential formulas in the language of exponential rings must define countable or co-countable sets in $\Cexp$.

\section{Zilber's exponential field \Bexp}

Zilber~\cite{Zilber05peACF0} constructed a QME exponential field, \Bexp, using the Hrushovski--Fra\"iss\'e amalgamation-with-predimension method together with the excellence technique. It is the unique model of cardinality continuum of some $L_{\omega,\omega_1}(Q)$-expressible axioms:
\begin{enumerate}
\item Algebraic properties of $\Cexp$:  $\mathrm{ACF}_0$, $\exp$ is a surjective homomorphism from the additive to the multiplicative group of the field, and the kernel is $\tau\Z$ for a transcendental $\tau$.
\item Schanuel's conjecture: for any tuple $\bar z$ we have $\td(\bar z, \exp(\bar z)) \geqslant \ldim_\Q(\bar z)$.
\item Strong Exponential-Algebraic Closedness (SEAC): For any algebraic subvariety $V \subseteq \ga^n \times \gm^n$ which is free and rotund, and defined over a finitely generated subfield $A$, there is $(\bar z,\exp(\bar z)) \in V$, generic in $V$ over $A$.
If we drop the genericity clause, this is called Exponential-Algebraic Closedness (EAC).
\item The natural pregeometry, exponential-algebraic closure, has the countable closure property: the closure of a finite set is countable. (This is true in $\Cexp$, proved using the Ax--Schanuel property.)
\end{enumerate}
For more details about the axioms and the logic needed to express them, see \cite{NAPE}.

A strengthening of Conjecture~\ref{QM conj} is:
\begin{conjecture}\label{strong QM conj}
$\Cexp \iso \Bexp$. More precisely, $\Cexp$ is a model of the axioms above. Equivalently, Schanuel's conjecture is true in $\Cexp$ and $\Cexp$ is  Strongly Exponentially-Algebraically Closed.
\end{conjecture}

D'Aquino, Macintyre and Terzo~\cite{DMT10,DMT14,DMT16} have done some analysis of what known properties of $\Cexp$ with analytic proofs can also be proved algebraically for \Bexp.

Although Schanuel's Conjecture seems well out of reach, proving EAC, or proving SEAC assuming Schanuel's conjecture to get the genericity of solutions, seems more plausible. Several people have made progress:

Marker~\cite{Marker06} proved the $n=1$ case. In unpublished work, Mantova and Masser proved the $n=2$ case.
Brownawell and Masser~\cite{BM17} proved the case of EAC where the projection of $V$ onto $\ga^n$ is dominant. See also \cite{DAFT18}. Another proof of this was given in \cite{AKM22} where a similar result for the exponential maps of abelian varieties was also given.

While working on the construction of $\Bexp$, Zilber was also trying to understand issues of uniformity needed to give first-order axioms, and this led him to his Conjecture on Intersections of subvarieties with Tori, or CIT, now known as the multiplicative part of the Zilber--Pink Conjecture~\cite{Zilber02esesc}. Conditional on this CIT, \cite{ECFCIT} gives an axiomatisation of the first-order theory of $\Bexp$.

In another direction, in 2008 I started a project with Martin Bays and Juan Diego Caycedo to do an analogous construction to $\Bexp$ but of a Weierstrass $\wp$-function. While Caycedo eventually left the project, this work led Martin and me to the papers previously mentioned on QME classes and eventually to \cite{PEM} where we rewrote and generalised the construction of $\Bexp$ to give our $\B_\wp$ and many other similar constructions. In that paper we also gave a strategy to prove Conjecture~\ref{QM conj} by showing:
\begin{theorem}\label{EAC implies QME}
If $\Cexp$ satisfies EAC then it is QME.
\end{theorem}

Using this method, and the density of the rationals in the reals, I was able to prove:
\begin{theorem}\cite{BCE} 
Let $\Gamma = \class{(z,\exp(z+q+2\pi i r))}{z \in \C, q,r \in \Q}$. Then the ``blurred exponential field'' $\tuple{\C;+,\cdot,\Gamma}$ is QME.
\end{theorem}

\section{Complex Powers}

As part of his PhD thesis~\cite{Gallinaro_thesis}, see also~\cite{Gallinaro22esetg}, Gallinaro recently proved that the part of EAC which relates to complex powers is true. Specifically, if $V$ is of the form $L \times W$ where $L \subseteq \ga^n$ is given by $\C$-linear equations and $W \subseteq \gm^n$ is any algebraic variety, and $L \times W$ is free and rotund, then there is $\bar z \in L(\C)$ such that $\exp(\bar z) \in W(\C)$.

The proof of Theorem~\ref{powers theorem} goes by using this EAC result together with the analogue of Theorem~\ref{EAC implies QME} for the powers setting, which is in \cite{GallinaroKirby23}.

In fact we prove QME for a slightly more expressive structure which also turns out to be easier to work with: the exponential sums language of Zilber from \cite{Zilber03powers,Zilber11tes}. Given a countable subfield $K$ of $\C$, he considers the two-sorted structure
\[ \C^K \quad \leteq \quad \left(\C_{K\text{-VS}}  \stackrel{\exp}{\longrightarrow} \C_{\text{field}} \right)\]
where the left sort has just the $K$-vector space structure and the right sort has the field structure. The induced structure on the right sort is that of $\tuple{\C;+,\cdot,(\Gamma_\lambda)_{\lambda \in K}}$, but not all automorphisms of that sort lift to the covering sort. The terms in the field sort with variables from the covering sort naturally correspond to complex exponential sums with exponents in $K$.

\newcommand{\etalchar}[1]{$^{#1}$}


\end{document}